\documentclass[12pt]{amsart}
\usepackage{graphicx}
\usepackage[active]{srcltx}
\newcounter{minutes}
\divide\time by 60
\newcounter{hours}
\multiply\time by 60 \addtocounter{minutes}{-\time}
 \makeatletter
\renewcommand*\subjclass[2][2000]{%
  \def\@subjclass{#2}%
  \@ifundefined{subjclassname@#1}{%
    \ClassWarning{\@classname}{Unknown edition (#1) of Mathematics
      Subject Classification; using '1991'.}%
  }{%
    \@xp\let\@xp\subjclassname\csname subjclassname@#1\endcsname
  }%
}
 \makeatother
\usepackage{enumerate,url,amssymb,  mathrsfs,yhmath}%, pdfsync}% ,showkeys, pdfsync}

\usepackage{enumerate,url,amssymb,  mathrsfs,yhmath}
\newtheorem{theorem}[equation]{\bf Theorem}

\newtheorem{lemma}[equation]{Lemma}
\newtheorem{proposition}[equation]{Proposition}
\newtheorem{example}[equation]{Example}

\newtheorem{remark}[equation]{Remark}

\numberwithin{equation}{section}

\def\XXint#1#2#3{{\setbox0=\hbox{$#1{#2#3}{\int}$}
\vcenter{\hbox{$#2#3$}}\kern-.5\wd0}}

\def\le{\leqslant}

\begin{document}

%\begin{center}
%\texttt{\small File:~\jobname .tex,
%          printed: \number\year-\number\month-\number\day,
%          \thehours.\ifnum\theminutes<10{0}\fi\theminutes}
%\end{center}

\title{On harmonic functions and the Schwarz lemma{\footnote{ File:~\jobname .tex,
          printed: \number\year-\number\month-\number\day,
          \thehours.\ifnum\theminutes<10{0}\fi\theminutes}}}

\author{David Kalaj}
\address{University of Montenegro, Faculty of Natural Sciences and
Mathematics, Cetinjski put b.b. 81000 Podgorica, Montenegro}
\email{davidk@t-com.me}

\author{Matti Vuorinen}
\address{Department of Mathematics, University of Turku, 20014
Turku, Finland} \email{vuorinen@utu.fi}

\begin{abstract}
We study the  Schwarz lemma for harmonic functions and prove sharp
versions for the cases of real harmonic functions and the norm of
harmonic mappings.
%We prove a version of Schwarz lemma for real harmonic functions and a version of Schwarz %lemma for the norm of harmonic mappings.
% The inequalities we prove are sharp and are extensions of some recent results.
\end{abstract}
\maketitle

%\tableofcontents

\section{Introduction}

We first recall that the hyperbolic metric $d_h(z,w)$ of the unit
disk $\mathbf{U}:=\{z:|z|<1\}$ is defined by
$$
\tanh \frac{d_h(z,w)}{2}= \frac{|z-w|}{|1-z \overline{w}|} \,.
$$
The classical Schwarz lemma says that an analytic function $f $ of
the unit disk into itself is a contraction in the hyperbolic metric,
i.e. for $z , w\in \mathbf{U}$
\begin{equation} \label{hypSLemma1}
\frac{|f(z)-f(w)|}{|1-f(z) \overline{f(w)}|} \le \frac{|z-w|}{|1-z
\overline{w}|}. \,
\end{equation}
Letting $z \to w$ we get
\begin{equation} \label{hypSLemma2}
|f'(z)| \le \frac{1-|f(z)|^2}{1-|z|^2} \,.
\end{equation}
In the particular case $f(0)=0$ we obtain from (\ref{hypSLemma1})
the following
\begin{lemma}{} \label{ClassSLemma}
  Let $f: {\mathbf U} \to {\mathbf U} $ be an analytic function
with $f(0)=0\,.$ Then $|f(z)| \le |z|$, $z \in {\mathbf U} \,.$
\end{lemma}

In this paper we will consider harmonic functions $f :{\mathbf U}
\to B_n:=\{x: |x|<1\}\subset \mathbf R^n$ i.e. functions whose
coordinate functions satisfy the Laplace equation. In this case, the
Cauchy-Riemann equations need not be satisfied. During the past
decade, harmonic functions have been extensively studied and many
results from the theory of analytic functions have been extended for
them. For instance, the Schwarz lemma has the following counterpart
in this context.

\begin{lemma}{{\rm (\cite{H}, \cite[p. 77]{D} and \cite[Theorem~4.4.6]{lib})}} \label{SLemma}
Let $f: {\mathbf U} \to {\mathbf U} $ be a harmonic functions with
$f(0)=0\,.$ Then $|f(z)| \le \frac{4}{\pi} {\mathrm{ tan}}^{-1}|z|$
and this inequality is sharp for each point $z \in {\mathbf U} \,.$
\end{lemma}

%Our first observation is that, although sharp, this lemma can be
%improved.
 Recently,
several authors have proved the following three extensions of the
Schwarz lemma, which is the motivation for our study.
\begin{proposition}\label{1}\cite[Theorem~6.26]{ABR}
Let $f$ be a real harmonic function of the unit disk into $(-1,1)$.
Then there hold the following sharp inequality
$$|\nabla f(0)|\le \frac{4}{\pi}.$$
\end{proposition}
\begin{proposition}\label{3}\cite[Proposition~3.1.]{MMM}
Let $f$ be a $K-$q.c. harmonic mapping of the unit disk into itself.
Then
$$|\nabla f(z)|\le K\frac{1-|f(z)|^2}{1-|z|^2}.$$
\end{proposition}
The previous proposition has its counterpart for the class of
quasiconformal hyperbolic harmonic mappings. See \cite{wan} for
related results.
\begin{proposition}\label{2}\cite[Theorem~1.1]{pav}
Let $f$ be an analytic function of the unit disk into the unit ball
$B_n\subset \mathbf C^n$. Then there holds the following sharp
inequality
$$|\nabla |f(z)||\le \frac{1-|f(z)|^2}{1-|z|^2}.$$
\end{proposition}
The results of this paper are:
\begin{theorem}\label{te1}
Let $f$ be a real harmonic function of the unit disk into $(-1,1)$.
Then  the following sharp inequality holds
\begin{equation}\label{sha}|\nabla f(z)|\le
\frac{4}{\pi}\frac{1-|f(z)|^2}{1-|z|^2},\ \ \ |z|<1.\end{equation}
\end{theorem}
\begin{theorem}\label{te22}
If $f$ is a harmonic mapping of the unit disk into the unit ball
$B_n\subset \mathbf R^n$, then
\begin{equation}\label{nor}\left|\nabla |f(z)|\right|\le
\frac{4}{\pi}\frac{1-|f(z)|^2}{1-|z|^2},\ \ \ |z|<1.\end{equation}
The inequality \eqref{nor} is sharp for each $z$.
\end{theorem}
\begin{theorem}\label{tete}
Let $f$ be a harmonic function of the unit disk into $(-1,1)$.  Then
we have $$d_h(f(z),f(w))\le \frac{4}{\pi}d_h(z,w),$$ where $d_h$ is
hyperbolic distance in  $(-1,1)$ and $\mathbf U$ respectively.
\end{theorem}
\begin{remark}
Theorem~\ref{te1} is an improvement and an extension of
Proposition~\ref{1}. Notice that
$$\frac{4}{\pi}\left({1-|f(0)|^2}\right)\le \frac{4}{\pi}.$$ Theorem~\ref{te22} is a
harmonic version of Proposition~\ref{2}.
\end{remark}
The following example is taken from \cite{kalaj.thesis} but it is
used there in some other context. It shows that we do not have any
complex version of Theorem~\ref{te1}, without imposing
quasiconformality (as is done in Proposition~\ref{3}).
\begin{example}
Let $z=re^{it}$ and let $f$ be a harmonic diffeomorphism of the unit
disk onto itself such that $\partial_r f$ exists in the boundary
point $1$, and $\partial_t f$ is not bounded on $z=1$. Then
$$\overline{\lim}_{r\to
1}\frac{1-|f(r)|^2}{1-r^2}=|f(1)|\left|\frac{\partial
f(r)}{dr}\right|_{r=1}=M<\infty$$ but
$$\sup_{r\in[0,1)}|\nabla f(r)|=\infty.$$
For instance: let  $\psi$ be the function defined on $[-\pi ,\pi]$
by
$$ \psi(\theta)=a\left(\theta+\int_0^{\theta}-\log|x|dx\right), $$
where $a$ is a positive constant such that
$$a\left(2\pi+\int_{-\pi}^{\pi}-\log|x|dx\right)=2\pi.$$
Take $F(\theta)=e^{i\psi(\theta)}$ and define $f(z)=P[F](z).$
\end{example}
\section{Proofs}
\begin{proof}[Proof of Theorem~\ref{te1}] Let  $h$ be the harmonic
conjugate of $f\,.$  Then $a=f+ih$ maps the unit disk $\mathbf U $
into the vertical strip $S= \{w: -1< \Re w< 1 \}\,.$
%The mapping
%$$ g(z) = \frac{2i}{\pi}\log \frac{1+z}{1-z}$$ maps the unit disk
%onto the vertical strip  Since $f$ is a real
%harmonic mapping of the unit disk into the interval $(-1,1)$, there
%exists an analytic function $a$ of the unit disk into the strip
%$-1\le \Re w\le 1$ such that $f(z)=\Re a(z)$.  Assume that $a(0)=0$
%i.e. $f(0)=0$.

Making use of the conformal map
$$ g(z) = \frac{2i}{\pi}\log \frac{1+z}{1-z}$$
of the unit disk $ \mathbf U$ onto the strip $S$ we see that
 there exists an analytic function
$b:  \mathbf U \to \mathbf U$, such that
$$a(z) =\frac{2i}{\pi}\log \frac{1+b(z)}{1-b(z)}.$$ By Schwarz
lemma for analytic functions we have
$$|b'(z)|\le \frac{1-|b(z)|^2}{1-|z|^2}.$$ On the other hand $$a'(z) =
\frac{4i}{\pi}\frac{b'(z)}{1-b(z)^2}$$ and $$|\nabla f|=|a'(z)|.$$
We will find the best possible constant $C$ such that
$$|\nabla f(z)|\le C \frac{1-|f(z)|^2}{1-|z|^2}.$$
As
$$|a'(z)|\le \frac{4}{\pi}\frac{1-|b(z)|^2}{|1-b(z)^2|} \, \frac{1}{1-|z|^2}$$
it will be enough to find the best possible constant $C$ such that
$$\frac{4}{\pi}\frac{1-|b(z)|^2}{|1-b(z)^2|} \, \frac{1}{1-|z|^2}\le C
\frac{1-|\Re a(z)|^2}{1-|z|^2}$$ or what is the same
$$\frac{4\pi}{({\pi^2-4|\mathrm{arg} \frac{1+b}{1-b}|^2})}\frac{1-|b|^2}{|1-b^2|}\le C\,, |b|<1\,.$$
Let $\omega= \frac{1+b}{1-b} = r e^{it}$. Then $-\pi/2\le t\le
\pi/2$ and
$$  1-|b|^2 = \frac{4 r  \cos  t}{r^2+ 2r \cos t +1},$$ $$
|1-b^2| = \frac{4 r |e^{it} |}{r^2+ 2r \cos t +1} ,$$ and hence the
last inequality with the constant $C= 4/\pi$ follows from
$$\frac{|\cos t|}{1-\frac{4}{\pi^2}t^2}\le 1$$
which holds for $t \in (-\pi/2,\pi/2) \, .$ This yields \eqref{sha}.

To show that the inequality \eqref{sha} is sharp, take the harmonic
function $$f(z) =
\frac{2}{\pi}\mathrm{arctan}\frac{2y}{1-x^2-y^2}.$$ It is easy to
see that $$|\nabla
f(0)|=\frac{4}{\pi}=\frac{4}{\pi}\frac{1-|f(0)|^2}{1-0^2}.$$
%Thus $C= \frac 4{\pi}$.
\end{proof}
\begin{proof}[Proof of Theorem~\ref{te22}]
Consider the function $g(z) = \left<f(z), b\right >$, where $|b|=1$,
$b\in\mathbf R^n$. Then by Theorem~\ref{te1}, for $|h|=1$, $h\in
\mathbf C$ we have
\begin{equation}\label{coin}\left<\nabla f(z) h, b\right>\le
\frac{4}{\pi}\frac{1-|g(z)|^2}{1-|z|^2}.\end{equation} Let
$$b=\frac{f(z)}{|f(z)|}$$ and $$h= \frac{\nabla f(z)^t b}{|\nabla f(z)^t
b|}.$$ Here $\nabla f(z)^t$ is the transpose of the matrix $\nabla
f(z)$. Then \eqref{coin} coincides with
$$|\nabla |f||\le \frac{4}{\pi}\frac{1-|f(z)|^2}{1-|z|^2}.$$ It
follows the desired conclusion.
\end{proof}
\begin{proof}[Proof of Theorem~\ref{tete}]
Let $m$ be a M\"obius transformation of the unit disk onto itself
that maps the points $0$ and $r>0$ to the points  $z$ and $w$. Then
the function $f\circ m$ is also harmonic. Denote it by $f$ as well.
Without loss of generality assume that $f(z)\le f(w)$. Let $C=\{t\in
[0,r]: f'(t)> 0\}$. Then $C = \bigcup_{n=0}^\infty C_n$, where $C_n$
are some intervals in $[0,r]$. Let $A_1=C_1$, $A_2=C_2\setminus\{t:
f(t)\in f(A_1)\}$, and by induction, let
$$A_{n+1} = \{t\in C_n: f(t) \not\in f(\bigcup_{k=1}^n A_k)\}.$$
Take $$A = \bigcup_{k=1}^\infty A_k.$$  Denote by $|\cdot |$ the
Lebesgue measure on real line. Then $[f(z),f(w)]\subset f([0,r])$
and $A$ is a measurable subset of $[0,r]$ such that $|[f(z),f(w)]|=
|f(A)|$ and $f$ is injective in $A$. Therefore by using \eqref{sha}
we have
\[\begin{split}d_h(f(z),f(w))&=\int_{f(z)}^{f(w)}\frac{|d\omega|}{1-|\omega|^2}\\&=\int_{A}\frac{|df(t)|}{1-|f(t)|^2}\\&\le
\frac{4}{\pi}\int_{A}\frac{|dt|}{1-|t|^2}\\&\le
\frac{4}{\pi}\int_{[0,r]}\frac{|dt|}{1-|t|^2}=
\frac{4}{\pi}d_h(z,w)\end{split}\]
\end{proof}
\subsection*{Acknowledgment} This work was completed during the visit of the
first author to the Universities of Helsinki and Turku in October
2010. The first author was partially supported by Ministry of
education of Montenegro. This research was supported by Academy of
Finland, project 2600066611.

\end{document}